\documentclass[reqno,12pt]{amsart}
\usepackage{txfonts}
\usepackage{mathrsfs}
\usepackage{amsmath}
\usepackage{amssymb}

\usepackage{amsthm}
\usepackage{verbatim}
\usepackage{latexsym, bm}

\title[Uniform Trace-free Bundles]{Uniform Trace-free Bundles of Triple Covers over Projective Planes}


\author[Xinyi Fang]{Xinyi Fang}

\address{School of Mathematical Sciences\\
Shanghai Key Laboratory of PMMP\\
East China Normal University\\
No. 500, Dongchuan Road\\
Shanghai, 200241, P. R. China}

\email{2315885681@qq.com}

\author[Yongtao Wang]{Yongtao Wang}

\address{School of Mathematical Sciences\\
Shanghai Key Laboratory of PMMP\\
East China Normal University\\
No. 500, Dongchuan Road\\
Shanghai, 200241, P. R. China}

\email{13564932708@163.com}

\thanks{The Research is Sponsored by National Natural Science Foundation of China (Grant No. 11531007, 11761141005) and Science and Technology Commission of Shanghai Municipality (Grant No. 18dz2271000).}

\theoremstyle{definition}
\newtheorem{theorem}[subsection]{Theorem}

\newtheorem{definition}[subsection]{Definition}
\newtheorem{proposition}[subsection]{Proposition}

\newtheorem{corollary}[subsection]{Corollary}
\newtheorem{remark}[subsection]{Remark}

\newfont{\drnew}{wncyr10}
\let\tilde=\widetilde

\allowdisplaybreaks

\def\dashfill{\leaders\hbox{\hbox to 3.25pt{\hrulefill}\hspace*{2pt}\hbox to 3.25pt{\hrulefill}}\hfill}
\newcommand{\CITE}[1]{{[#1]}}
\let\cite=\CITE

\renewcommand{\O}{\mathcal{O}}
\renewcommand{\P}{\mathbb{P}}

\newcommand{\DP}{\Delta_\pi}
\newcommand{\OP}{\Omega_{\P^2}}
\newcommand{\TP}{T_{\P^2}}

\begin{document}

%


\maketitle
\begin{abstract}
Studying coverings over algebraic varieties is an effective method in algebraic geometry. By combining the technique of triple cover from Miranda and Tan, we proved that if the degree of the branch divisor of a normal triple cover over $\P^2$ is no more than $18$ and not equal to $16$, then the trace-free bundle of the triple cover is uniform. Moreover, we totally classified all these trace-free bundles.

\end{abstract}

\section{Introduction}
It is classically known that every holomorphic vector bundle on the projective line splits as a direct sum of line bundles. However, for $n$-dimensional projective spaces ($n\ge 2$), the situation is very complicated. So the splitting of vector bundles on higher dimensional projective spaces has long been a major concerning among the problems on vector bundles in algebraic geometry. The typical example of indecomposable $2$-bundle over $\P^2$ is the tangent bundle $\TP$ or its dual. For $n\ge 3$, unsplit $2$-bundles over complex projective space of high-dimension are difficult to construct. Tango (\cite{Tango}) constructed holomorphic $(n-1)$-bundle over  $\P^n$ for any $n\ge 3$. In 1973, Horrocks and Mumford found an unsplit $2$-bundle over $\P^4$ (\cite{HM}). The bundle is essentially the only known non-splitting holomorphic $2$-bundle over $\P^n(n\ge4)$. In 1974, Hartshorne conjectured that for $n\ge7$, every holomorphic $2$-bundle over complex projective space $\P^n$ splits. However, this conjecture still remains open. For non-splitting vector bundles, we may obtain partial classification results after restricting to certain classes of vector bundles. One of the classes that has been studied more extensively is \emph{uniform} vector bundles, namely those in which the splitting type is independent of the chosen line. The notion of the uniform vector bundle appears first in a paper of Schwarzenberger (\cite{Schw}) and lots of works are about classification of such bundles.

Triple covers are effective tools for constructing rank $2$ vector bundles over $\P^n$ because the trace-free sheaf is of rank $2$ naturally. Triple covers were first systematically studied by Miranda (\cite{M}). Later, many mathematicians performed more detailed work (\cite {TF,TS2}). By studying the integral closure of a cubic extension, Tan gave a more detailed triple cover data such that the structure of triple covers become more clear.

 \begin{definition}
 Let $ \pi : X \rightarrow \P^2 $ be a normal triple cover. Denote the branch locus of the triple cover over $\P^2$ by $\DP$, which is a subvariety of pure dimension $1$ over $\P^2$. Regard it as a reduced divisor over $\P^2$.  We may write $\DP=S_\pi + T_\pi $, where $S_\pi$ (resp. $T_\pi$) is of ramification index $2$ (resp. $3$).  Denote $S_\pi + 2 T_\pi $ by $ \overline{\DP} $, which is called the \textit{branch divisor} of $\pi$.
 \end{definition}

 Over the projective plane, the trace-free sheaf of a triple cover is determined by the branch locus.
        It is an interesting problem to determine the trace-free sheaf of a triple cover by the degree of branch divisor.
        Moreover, we can give the geometry of the given triple covers.

                If $ S_\pi = 0 $, then the normal triple cover $ \pi : X \rightarrow \P^2 $ is Galois.
        Miranda has proven that in this case, the trace-free sheaf splits.
        More precisely, a triple cover is Galois if and only if $ a = d = 0 $ \cite{M} (The symbol will be introduced in chapter \ref{2}).
        For general situation, Tokunaga solved the case of $(\deg S_\pi, \deg T_\pi) = (2,1), (2,2), (4,0), (4,1)$ by the theory of dihedral cover \cite{HT1,HT2}. In 2010, Taketo proved that if the degree of branch divisor is $6$, then the trace-free sheaf is uniform and the surface is either the cubic surface in $\P^3$, or the pull back of the dual curve correspond to the discriminant (\cite{TS}).

In this paper, by combining the technique of triple cover from Miranda and Tan, we proved that if the degree of the branch divisor of a normal triple cover over $\P^2$ is no more than $18$ and not equal to $16$, then the trace-free bundle of the triple cover is uniform.

 \begin{theorem}\label{main}

            Let $ \pi : X \rightarrow \P^2 $ be a normal triple cover, $ \deg \overline{\DP} \leqslant 18 $, and $ \deg \overline{\DP} \not = 16 $. Then $E_\pi$ is uniform. Moreover,

            \begin{description}

                  \item[$ \deg \overline{\DP} = 2 $],
                  $ E_\pi \cong \OP(1) $;

                  \item[$ \deg \overline{\DP} = 4 $],
                  $ E_\pi \cong \O_{\P^2}(-1) \oplus \O_{\P^2}(-1) $;

                  \item[$ \deg \overline{\DP} = 6 $],
                  $ E_\pi \cong \O_{\P^2}(-1) \oplus \O_{\P^2}(-2) $ or $ E_\pi \cong \OP $;

                  \item[$ \deg \overline{\DP} = 8 $],
                  $ E_\pi \cong \O_{\P^2}(-2) \oplus \O_{\P^2}(-2) $;

                  \item[$ \deg \overline{\DP} = 10 $],
                  $ E_\pi \cong \O_{\P^2}(-2) \oplus \O_{\P^2}(-3) $ or $ E_\pi \cong \OP(-1) $;

                  \item[$ \deg \overline{\DP} = 12 $],
                  $ E_\pi \cong \O_{\P^2}(-3) \oplus \O_{\P^2}(-3) $ or $ E_\pi \cong \O_{\P^2}(-2) \oplus \O_{\P^2}(-4) $;

                  \item[$ \deg \overline{\DP} = 14 $],
                  $ E_\pi \cong \O_{\P^2}(-3) \oplus \O_{\P^2}(-4) $ or $ E_\pi \cong \OP(-2) $;

                  \item[$ \deg \overline{\DP} = 18 $],
                  $ E_\pi \cong \O_{\P^2}(-4) \oplus \O_{\P^2}(-5) $ or $ E_\pi \cong \O_{\P^2}(-3) \oplus \O_{\P^2}(-6) $ or $ E_\pi \cong \OP(-3) $.

            \end{description}

        \end{theorem}

\section{Preliminaires}\label{2}
 \begin{definition}

            We call a finite flat morphism $ \pi : X \rightarrow Y $ from a scheme $X$ to a variety $Y$ of degree three a triple cover over $Y$.
            If in addition, $X$ and $Y$ are normal varities, we call it a normal triple cover over $Y$.
        \end{definition}

        In this paper, we mainly talk about normal triple covers over $\P^2$.

        \begin{definition}

            Let $ \pi : X \rightarrow \P^2 $ be a normal triple cover.
            The kernel of the trace map $ \pi_*\O_X \rightarrow \O_{\P^2} $ is a rank $2$ locally free $\O_{\P^2}$-module, denoted by $E_\pi$, called the \textit{trace-free sheaf} of $ \pi : X \rightarrow \P^2 $.
            Particularly, $ \pi_*\O_X = \O_{\P^2} \oplus E_\pi $.

        \end{definition}

 \begin{remark}
Trace-free sheaves of triple covers over projective planes are precisely vector bundles, for trace-free sheaves are reflexive sheaves over projective planes and the singularity set of a reflexive sheaf is of codimension at least 3.
\end{remark}

        \begin{theorem} \label{2.2}

            (\cite{M}, Theorem 3.6)
            For every given $2$-bundle $\mathcal{E}$ over $\P^2$, the triple cover, whose trace-free sheaf $ E_\pi = \mathcal{E} $, given by the $\O_{\P^2}$ - algebra structure of $ \mathcal{A} = \O_{\P^2} \oplus \mathcal{E} $ is one-to-one corresponding to a $\O_{\P^2}$-linear homomorphism $ \Phi : S^3 \mathcal{E} \rightarrow \det \mathcal{E} $.

        \end{theorem}

        We can discribe the correspondence more precisely.
        Choose $\{ z,w \}$ as a group of local basis of $\mathcal{E}$. Denote the homomorphism given by the multiplication of  $\mathcal{A}$ by $ \Phi : S^2\mathcal{E} \rightarrow \mathcal{A} $.
        Then $\phi$ has the following form:

        \begin{align*}
            \phi(z^2) &= az + bw + 2e, \\
            \phi(zw) &= - dz - aw -f, \\
            \phi(w^2) &= cz + dw + 2g.
        \end{align*}

        where $a,b,c,d \in \O_{\P^2}$, $e = a^2 - bd,\ f = ad - bc,\ g = d^2 - ac $.
        Particularly, $\mathcal{A}$ is an intergral domain, $b \not =0$ and $c \not =0$.

        Define $\Phi : S^3\mathcal{E} \rightarrow \det \mathcal{E}$. On the basis, given by
        $$ \Phi(z^3) = -b(z \wedge w), \quad \Phi(z^2w) = a(z \wedge w), $$
        $$ \Phi(zw^2) = -d(z \wedge w), \quad \Phi(w^3) = c(z \wedge w). $$

        It is not difficult to verify that it is irrelevant with the selection of $\{z,w\}$.

        \begin{proposition} \label{tt}

            (\cite{M}, Lemma 4.5, Proposition 4.7; \cite{TSL}, Theorem 1.3)
            Let $ \pi: X \rightarrow \P^2 $ be the triple cover given by the homomophism $ \Phi:S^3\mathcal{E} \rightarrow \det \mathcal{E} $ in the above way.
            Then the branch divisor $\overline{\DP}$ is locally given by
            $$ D : = f^2 - 4eg=0, $$
            where the definition of $ e, f, g $ is the same as above. In addition, the line bundle corresponds to the $\overline{\DP}$ is $ (\det E_{\pi})^{-2} $.

        \end{proposition}

        \begin{proposition}\label{mm}

            (\cite{M}, Section 6)
            Let $ \pi : X \rightarrow \P^2 $ be a normal triple cover.
            If $ E_\pi \cong \O_{\P^2}(t_1) \oplus \O_{\P^2}(t_2) $, then
            $ a \in H^0(\O_{\P^2}(-t_1)), b \in H^0(\O_{\P^2}(-2t_1 + t_2)), c \in H^0(\O_{\P^2}(t_1 - 2t_2)), d \in H^0(\O_{\P^2}(-t_2)) $.
            Therefore $ 2t_1 \leqslant t_2, 2t_2 \leqslant t_1$.

        \end{proposition}

        \begin{proposition} \label{2.3}

            (\cite{TF})
            Let $ \pi : X \rightarrow \P^2 $ be a triple cover, $\pi$ not $\mathrm{\acute{e}tale}$. Then $X$ is a cubic section of the total space of a line bundle $\mathcal{L}$ over $\P^2$ if and only if $E_\pi \cong \mathcal{L}^{-1} \oplus \mathcal{L}^{-2}$.

        \end{proposition}

 By studying the integral closure of a cubic extension, Tan gave a more detailed triple cover data such that the structure of triple covers become more explicit and clear (see \cite{TS2}).

        \begin{definition}

 Let $\mathcal{L}$ be a line bundle on $\P^2$ and the minimal equation of a normal triple cover $ \pi : X \rightarrow \P^2 $ is
            $$ z^3 + sz + t = 0$$
            where $ 0 \not = s \in H^0(\P^2, \mathcal{L}^2) $, $ 0 \not = t \in H^0(\P^2, \mathcal{L}^3) $.
            Otherwise the triple cover is degenerated or cyclic.
            Because of the minimality of the equation, there is a unique decomposition:
            $$ s = a_1 a_2^2 b_1 a_0 , \quad t = a_1 a_2^2 b_1^2 b_0 $$
            where $ a_1, a_2, b_1 $ is square-free, $ (a_1, a_2) = 1 $, $ (a_i, b_j) = 1 $, $ i = 0, 1, 2, j = 0, 1 $. Denote
            $$ a = 4 a_1 a_2^2 a_0^3, \quad b = 27 b_1 b_0^2, \quad c = a + b = c_1 c_0^2 $$
            where $c_1$ square-free.
            Then
            $$ a = \frac{4 s^3}{\gcd (s^3,t^2)}, \quad b = \frac{27 t^2}{\gcd (s^3,t^2)}, \quad c = \frac{4 s^3 + 27 t^2}{\gcd (s^3,t^2)}. $$

            Call $ (a, b, c) $ the $ abc $ data of the triple cover $ \pi $.

        \end{definition}

        Denote
        $$ g_1 = \frac{2}{3} a_0 a_2, \quad g_2 = b_0, \quad g_3 = c_0. $$

        Use capital letters to denote the divisors of the corresponding sections.

        \begin{theorem}

            (\cite{TSL})
            Let $ \pi : X \rightarrow \P^2 $ be a normal triple cover,

            \begin{enumerate}
                \item the divisor over which $\pi$ is totally ramified is $ A_1 + A_2 $.
                    The divisor over which $\pi$ is simple ramified is  $ B_1 + C_1 $.
                    The branch locus of $\pi$ is $ 2A_1 + 2A_2 + B_1 + C_1 $.
                \item the image in $X$ of non-normal locus is $ A_2 + B_1 + C_0 $.
            \end{enumerate}

        \end{theorem}

        We will discribe the transformation between two types of triple cover data.  We use $\tilde{\cdot}$ to denote the section or divisor in Miranda triple cover data. E.g. $\tilde{a}$, $\tilde{b}$, $\tilde{A}$, $\tilde{B}$.

            By direct computation, we can get the following corollary.

            \begin{corollary} \label{2.6}
Denote $ U_1 = \P^2 \backslash C_0 $ and $ U_2 = \P^2 \backslash (A_0 + A_2) $.
                Let $ \{ \tilde{a},\tilde{b},\tilde{c},\tilde{d} \} $ be the Miranda triple cover data, $ \{ a_0, a_1, a_2, b_0, b_1 \} $ be the S. -L. Tan triple cover data, we can get the transformation between two types of data:

                Over $U_1$,
                $$ \tilde{a} = 0, \quad \tilde{b} = a_2 b_1, \quad \tilde{c} = -a_1 b_0, \quad \tilde{d} = \frac{a_0 a_1 a_2}{3}. $$

                Over $U_2$,
                $$ \tilde{a} = \frac{3 b_0 b_1}{a_0}, \quad \tilde{b} = \frac{3 b_1 c_0}{2 a_0}, \quad \tilde{c} = \frac{b_0 c_1}{18 a_0}, \quad \tilde{d} = \frac{c_0 c_1}{18 a_0}.  $$

            \end{corollary}


     Let $E$ be a rank $2$ vector bundle over $\P^2$. According to the theorem of Grothendieck for every line $l\in {\P^2}^{*}$, there is a $2$-tuple \[a_E(l)=(a_1(l),a_2(l))\in\mathbb{Z}^{2};~ a_1(l)\geq a_2(l)\]
     with $E|L\cong \O_{L}(a_1(l)) \oplus \O_{L}(a_2(l)) $. We give $\mathbb{Z}^{2}$ the lexicographical ordering. Let\[
     \underline{a}_E=\inf_{l\in {\P^2}^{*}}a_E(l)
     \]
     \begin{definition}
     	$\underline{a}_E$ is the generic splitting type of $E$, $S_E=\{l\in {\P^2}^{*}|a_E(l)>\underline{a}_E\}$ is the set of jump lines.
     \end{definition}
     \begin{remark}
     	$U_E:={\P^2}^{*}\backslash S_E$ is a non-empty Zariski-open subset of ${\P^2}^{*}$.
     \end{remark}
     \begin{definition}
     	$E$ is uniform if there is no jumping line.
     \end{definition}
     For uniform $2$-bundles over $\P^2$, there is a specific depiction.
     \begin{theorem} (Van de Ven) \label{2.1}
     	
     	Every uniform $2$-bundle over $\P^2$ either splits, or is isomorphic to $\OP(t)(t\in \mathbb{Z} )$.
     	
     \end{theorem}

\section{Proof of the Main Theorem}

    \subsection{The proof of theorem \ref{main}}

        \begin{proof}

            Denote by $(E_\pi,\Phi)$ the triple cover data of $\pi$.

            As $ \deg \overline{\DP} = \deg S_\pi + 2 \deg T_\pi $ is even, for $S_\pi$ is a even divisor. Let $ \deg \overline{\DP} = 2k $, then by proposition \ref{tt} and hypothesis,
\[\det E_\pi \cong \O_{\P^2}(-k), ~k \in \{ 1,2,3,4,5,6,7,9 \} .\]

            {\bf{Case 1}} \quad $k$ is even (i.e. $ k = 2,4,6 $)

            Assume $ k = 2t $, where $ t \in \{ 1,2,3 \}$.
            Suppose $E_\pi$ is not uniform and $L$ is a jump line, $L$ is defined by $u_1 = 0$.
            There is some positive integer $m$ such that $ E_\pi|L \cong \O_L(-t + m) \oplus \O_L(-t - m) $.
            $\Phi|_L$ gives the section due to proposition \ref{mm}:
            $$ \tilde{a}|_L \in H^0(\O_L(t - m)), \quad \tilde{b}|_L \in H^0(\O_L(t - 3m)), $$
            $$ \tilde{c}|_L \in H^0(\O_L(t + 3m)), \quad \tilde{d}|_L \in H^0(\O_L(t + m)). $$

            Notice $ \tilde{b}|_L = 0 $ or $ \tilde{b}|_L $ is non-zero constant.

            If $ \tilde{b}|_L = 0 $ , restrict $\pi$ to the line $L$.
            By corollary \ref{2.6}, we have $ a_2|_L=0$ or $c_0|_L = 0 $ or $ b_1|_L = 0 $, i.e. $ u_1 \mid a_2$ or $u_1 \mid c_0 $ or $ u_1 \mid b_1 $.
 \begin{description}

                  \item[1. ] $ u_1 \mid a_2$ or $u_1 \mid c_0 $

                  Consider the $abc$ data, we always have $\tilde{d}|_L=0 $. This contradicts the previous statement $\tilde{d}|_L \in H^0(\O_L(t + m))$.

                  \item[2. ] $ u_1 \mid b_1 $

                  As $ b_1|_L = 0 $, i.e. the equation of triple cover over $L$ is $z^3|L= 0 $. In this case, $\pi$ is totally ramified over $L$.
                  The totally ramified of $\pi$ is $ A_1 + A_2 $, that is $ u_1 \mid a_1 a_2 $, which contradict with the coprime condition.

            \end{description}

            If $ \tilde{b}|_L $ is non-zero constant (the degree of the branch divisor is $12$), then $ \tilde{b} $ is constant over the whole $\P^2$. Due to corollary \ref{2.6}, it is easy to get $ a_0, a_2, b_1, c_0 $ are all constant.
             Therefore the $abc$ equation is given by 
             \begin{equation} \label{eq1}
             a_1 + b_0^2 = c_1.
             \end{equation}
          Since  $ \tilde{b}|_L $ is non-zero constant, the generic splitting type of $E_\pi$ over general line is $(-3,-3)$.
            Consider the restriction of triple cover to a general line $L'$.
            By comparing the degree of equation (\ref{eq1}) restricting to $L'$, we find that the degree of the left side is $\deg b_0^2|_{L'} = 2 \deg b_0|_{L'} $  while the degree of the right side is $ \deg c_1|_{L'} = 3 $, which is a contradiction.

            In summary, when $k$ is even, the trace-free sheaf $E_\pi$ is uniform.

            {\bf{Case 2}} \quad $k$ is odd (i.e. $ k = 1,3,5,7,9 $)

            Assume $ k = 2t - 1 $, where $ t \in \{ 1,2,3,4,5 \}$.
            Suppose $E_\pi$ is not uniform and $L$ is a jump line, $L$ is defined by $u_1 = 0$.
            There is some positive integer $m$ such that $ E_\pi|L \cong \O_L(-t + m+1) \oplus \O_L(-t - m) $.
            $\Phi|_L$ gives the section due to proposition \ref{mm}:
            $$ \tilde{a}|_L \in H^0(\O_L(t - 1 - m)), \quad \tilde{b}|_L \in H^0(\O_L(t - 2 - 3m)), $$
            $$ \tilde{c}|_L \in H^0(\O_L(t + 3m + 1)), \quad \tilde{d}|_L \in H^0(\O_L(t + m)). $$

            Notice that $ \tilde{b}|_L = 0 $ or $ \tilde{b}|_L $ is non-zero constant.
            The rest of the proof is similar.

            In summary, when the condition is met, $E_\pi$ is uniform.

            Thanks for the result of Van de Ven (theorem \ref{2.1}), the structure of rank $2$ uniform bundle is totally known. By computating the first chern class of the trace-free vector bundle, we get

            \begin{align*}
                \deg \overline{\DP} & = 2  \quad E_\pi \cong \OP(1), \\
                \deg \overline{\DP} & = 4  \quad E_\pi \cong \O_{\P^2}(-1) \oplus \O_{\P^2}(-1), \\
                \deg \overline{\DP} & = 6  \quad E_\pi \cong \O_{\P^2}(-1) \oplus \O_{\P^2}(-2) \ or \ E_\pi \cong \OP, \\
                \deg \overline{\DP} & = 8  \quad E_\pi \cong \O_{\P^2}(-2) \oplus \O_{\P^2}(-2), \\
                \deg \overline{\DP} & = 10 \quad E_\pi \cong \O_{\P^2}(-2) \oplus \O_{\P^2}(-3) \ or \ E_\pi \cong \OP(-1), \\
                \deg \overline{\DP} & = 12 \quad E_\pi \cong \O_{\P^2}(-3) \oplus \O_{\P^2}(-3) \ or \ E_\pi \cong \O_{\P^2}(-2) \oplus \O_{\P^2}(-4), \\
                \deg \overline{\DP} & = 14 \quad E_\pi \cong \O_{\P^2}(-3) \oplus \O_{\P^2}(-4) \ or \ E_\pi \cong \OP(-2), \\
                \deg \overline{\DP} & = 18 \quad E_\pi \cong \O_{\P^2}(-4) \oplus \O_{\P^2}(-5) \\
                & or \ E_\pi \cong \O_{\P^2}(-3) \oplus \O_{\P^2}(-6)  \  or \ E_\pi \cong \OP(-3).
            \end{align*}

        \end{proof}

%
%
%
%
%
%
%
        For most of the above cases, there are some precise descriptions of the covering spaces by Miranda's results (See \cite{M}).
        
        \begin{corollary}
            Let $ \pi : X \rightarrow \P^2 $ be a normal triple cover, $ E_\pi $ be the corresponding trace-free bundle.
\begin{enumerate}
            \item $\deg \overline{\DP} = 4$, the surface $X$ is the Steiner cubic in $\P^4$ and the triple cover map is projection;

            \item $\deg \overline{\DP}= 6$ and $E_\pi$ splits, $X$ is a cubic hypersurface in $\P^3$. Otherwise $X$ is birationally equivalent to $B\times \P^1$, where $B$ is a elliptic curve;

            \item $\deg \overline{\DP}= 8$, the surface $X$ is the blow-up of $\P^1 \times \P^1$ at nine points; the map to the plane is given by curves of bidegree (2,3) through the nine base points;

            \item $\deg \overline{\DP}= 10$ and $E_\pi$ splits, $X$ is a quartic surface blown up at one point and the map is projection from the point. Otherwise $X$ is the $13$-fold blow-up of the plane,mapped to $\P^2$ via quartics through $13$ base points which impose only $12$ conditions on quartics;

           \item $\deg \overline{\DP}= 12$ and $ E_\pi \cong \O_{\P^2}(-2) \oplus \O_{\P^2}(-4) $, $X$ is a surface of general type with $ p_g(X) = K_X^2 = 3 $ and $f$ is the canonical map. Otherwise $X$ is an elliptic surface over $\P^1$ (the elliptic structure being given by the canonical map) and the triple covering is defined by a linear system of genus 4 trisections of the elliptic structure;

            \item $\deg \overline{\DP}=14$ and $E_\pi$ splits, $X$ is a quintic surface in $\P^3$ with a double point \textit{p} and \textit{f} is projection from \textit{p}. Otherwise $X$ is a surface with $q(X)=0$, $p_g(X)=3$, $K_X^2=2$ and $e(X)=46$.
\end{enumerate}
        \end{corollary}



\begin{thebibliography}{99}

    \thispagestyle{plain}

    \bibitem[TF]{TF}
    \newblock T. Fujita,
    \newblock Triple covers by smooth manifolds,
    \newblock \emph{J. Fac. Sci. Univ. Tokyo}, Sect. IA Math 35(1988), 169--175.

    \bibitem[HM]{HM}
    \newblock G. Horrocks, D. Mumford,
    \newblock A rank 2 vector bundle on $\P^4$ with 15000 symmetries,
    \newblock \emph{Topology}, 12(1973), 63--81.

    \bibitem[M]{M}
    \newblock R. Miranda,
    \newblock Triple covers in algebraic geometry,
    \newblock \emph{Amer. J. Math.}, 107(1983), 1123--1158.

    \bibitem[VB]{VB}
    \newblock Chirstian Okonek, Michael Schneider and Heinz Spindler,
    \newblock \emph{Vector Bundles on Complex Projective Spaces},
    \newblock New York: Springer, 2011.

    \bibitem[Schw]{Schw}
    \newblock R. L. E. Schwarzenberger,
    \newblock \emph{Vector bundles on the projective plane},
    \newblock  \emph{Proc. London Math. Soc.}, 11(3)(1961), 623–640.

    \bibitem[TS]{TS}
    \newblock T. Shirane,
    \newblock A note on normal triple covers over $\P^2$ with branch divisors of degree 6,
    \newblock \emph{Kodai Mathematical Journal}, 107(2010), 1123--1158.

    \bibitem[TS2]{TS2}
    \newblock T. Shirane,
    \newblock Families of Galois closure curves for plane quintic curves,
    \newblock \emph{J. Algebra}, 342(2011), 175--196.

    \bibitem[TSL]{TSL}
    \newblock Sheng-Li Tan,
    \newblock Triple covers on smooth algebraic varieties,
    \newblock \emph{Stud. Adv. Math.}, Vol. 29(2002), 143--164.

    \bibitem[TSL2]{TSL2}
    \newblock Sheng-Li Tan,
    \newblock Intergral closure of a cubic extension and applications,
    \newblock \emph{Proc. Amer. Math. Soc.}, 129(2001), no.9, 2553--2562.


    \bibitem[Tango]{Tango}
    \newblock H. Tango,
    \newblock An example of indecomposable vector bundle of rank $n−1$ on $\P^n$,
    \newblock \emph{J. Math. Kyoto Univ.},16(1976), 137–141.

    \bibitem[HT1]{HT1}
    \newblock H. Tokunaga,
    \newblock On dihedral Galois covering,
    \newblock \emph{Canad. J. Math.},46(1994), 1299--1317.

    \bibitem[HT2]{HT2}
    \newblock H. Tokunaga,
    \newblock Dihedral covers and an elementary arithmetic on elliptic surfaces,
    \newblock \emph{J. Math. Kyoto Univ},44(2004), 255--270.

    \bibitem[AV]{AV}
    \newblock A. Van de Ven,
    \newblock On uniform vector bundles,
    \newblock \emph{Math. Ann.}, 195(1972), 245--248.


\end{thebibliography}
\end{document}